\documentclass{amsart}

\usepackage{amsmath,amsfonts,amssymb,hyperref,mdwlist,graphicx,enumitem,color,multicol}
\usepackage[cmtip,all]{xy}
\usepackage[utf8]{inputenc}
\usepackage{cite}
\newtheorem{theorem}{Theorem}[section]
\newtheorem{lemma}[theorem]{Lemma}
\newtheorem{proposition}[theorem]{Proposition}
\theoremstyle{definition}
\newtheorem{definition}[theorem]{Definition}

\theoremstyle{remark}

\numberwithin{equation}{section}

\newcommand{\vect}{\textnormal{\bf{span}}\,} 
\newcommand{\corf}{\mathbb{F}}\newcommand{\corfd}{\corf\cdot} 
\newcommand{\rad}{\mathcal{N}}\newcommand{\radd}{\mathcal{N}\,^2}
\newcommand{\bim}{\mathcal{M}}\newcommand{\obim}{\mathcal{M}^{\textnormal{op}}}

 \newcommand{\reg}{\mathcal{R}\textnormal{eg}\,} 
\newcommand{\am}{^{(+)}}\newcommand{\wt}{\widetilde}
\newcommand{\alga}{\mathcal{A}}
\newcommand{\algs}{\mathcal{S}}\newcommand{\algb}{\mathcal{B}}
\newcommand{\jor}{\mathfrak{J}}

\newcommand{\clasekmj}{\mathcal{K}(\jor,\mathfrak{M})} \newcommand{\mj}{\mathfrak{M}(\jor)}
\newcommand{\heraa}{\mathcal{H}(\mathcal{A},\ast)} 

\newcommand{\alMnn}{\mathcal{M}_{n|n}(\corf)}

\newcommand{\osp}{\textnormal{osp}}\newcommand{\trp}{\textnormal{trp}} 
\newcommand{\suMnm}{\mathcal{M}_{n|m}(\corf)\am}
\newcommand{\josnm}{\jor\osp_{n|2m}(\corf)}
\newcommand{\jp}{JP_n}
\newcommand{\jpd}{JP_2}

\newcommand{\algdt}{\mathcal{D}_{t}}
\newcommand{\kap}{\mathcal{K}_3} 
\newcommand{\kac}{\mathcal{K}_{10}}\newcommand{\spformaf}{\jor(V,f)} 
\newcommand{\eun}{e^{1n}}\newcommand{\enu}{e^{n1}}\newcommand{\enn}{e^{n}}
\newcommand{\wh}{\widetilde{h}} \newcommand{\bh}{\bar{h}}\newcommand{\wu}{\widetilde{u}}
\newcommand{\ws}{\widetilde{s}} \newcommand{\bs}{\bar{s}}
 
\newcommand{\hh}{\widehat{h}} 
\newcommand{\hs}{\widehat{s}} 

\begin{document}
\title[Jordan superalgebra of type $\jp$, $n\geq 3$ and the WPT]
{Jordan superalgebra of type $\jp$, $n\geq 3$ and the Wedderburn principal theorem}
\author{G\'omez-Gonz\'alez F. A. } 
\address{Instituto de Matemáticas, Universidad de Antioquia, Colombia}
\email{faber.gomez@udea.edu.co}
\author{Ram\'irez Berm\'udez J.A. } 
\address{Instituto de Matemáticas, Universidad de Antioquia, Colombia} 
\email{jalexander.ramirez@udea.edu.co} 
\thanks{The first author was partially supported by Universidad de Antioquia, CODI 2016-12949.}

\begin{abstract}
We investigate an analogue to the Wedderburn Principal Theorem (WPT) for 
a finite-dimensional Jordan superalgebra 
$\jor$ with solvable radical 
$\rad$ such that 
$\radd=0$ and $\jor/\rad\cong\jp$, $n\geq 3$.
 
We consider $\rad$ as an irreducible $\jp$-bimodule 
and we prove that the WPT holds for $\jor$.
\end{abstract}
\maketitle
\noindent {\bfseries{Key Words:}} Decomposition theorem; Jordan superalgebra; Wedderburn Principal Theorem, split null extension.

\medskip\noindent
{\bfseries{2010 Mathematics Subject Classification:}} 17A70; 17A15; 17C50.

\section{Introduction}
It is known that for every finite-dimensional associative algebra $\alga$ 
over a field $\corf$ of characteristic zero, 
with radical $\rad$ there exists a subalgebra $\algs$
such that $\alga=\algs\oplus\rad$ and $\algs\cong\alga/\rad$.
For the case of algebras over the field of complex numbers, 
this result was proved by T. Molien \cite{Mol,Mol2} and 
generalized over an arbitrary field of characteristic zero by 
S. Epsteen and J. H. Maclagan-Wedderburn \cite{Wed1}.
This result is known as the \emph{Wedderburn Principal Theorem} (WPT). 
For Jordan algebras, the WPT was proved  by 
A. A. Albert \cite{Alb1}, A. J. Penico \cite{Pen}, and 
V. G. Azkinuze \cite{Azk}. In the case of alternative algebras, the WPT was proved by 
R. D. Schafer \cite{Sch}.

Alternative superalgebras were introduced by 
E. Zelmanov and I. P. Shestakov \cite{ZelShesalternativas} 
and Jordan superalgebras were introduced by 
V. G. Kac \cite{Kac} and I. Kaplansky \cite{Kap1}.
One may ask whether  an analogue to the WPT is still true for 
finite-dimensional superalgebras. 
For alternative superalgebras over a field of characteristic zero,  
N. A. Pisarenko \cite{Pis} proved that the WPT holds 
when some conditions are imposed on the types of simple summands in semi-simple quotient. 
Besides, N. A. Pisarenko proved that the conditions imposed are essential, 
these were illustrated by counterexamples (see \cite{Pis}). 
In the case of alternative superalgebras over a field of characteristic 3, 
some counterexamples to the WPT were constructed by 
M. C. Lopez Diaz \cite{LopD}.
In the case of Jordan superalgebras with solvable radical 
$\rad$ such that $\radd=0$, the first author  \cite{Fgg18} 
showed that it is possible to reduce the problem  to 
 simple quotients of Jordan superalgebras. 
Moreover, he proved that if radical $\rad$
satisfies $\rad^2=0$ 
then it is enough to verify only the cases when $\rad$ 
is an irreducible Jordan $\jor$-superbimodule.
Then, the proof of the WPT can be completed considering simple Jordan superalgebras 
 $\jor$ case by case. An analogue to the WPT holds
only when some conditions are imposed on the irreducible 
 $\jor$-superbimodules contained in the radical $\rad$. 
 It  has been shown
 that the restrictions are essential (see \cite{Fgg18,Fgg16,FggVr}).
 
Using the notation given by C. Martinez and E. I. Zelmanov  in \cite{zelmar2}, 
the first author  \cite{Fgg18, Fgg16} studied the WPT for the cases 
when semisimple part of $\jor$ is isomorphic to
(sometimes we say that $\jor$ is of type)
 $\spformaf$, $\algdt$, $\kac$,  $\kap$,
 $\suMnm$ and $\mathcal{M}_{1\mid 1}(\corf)^{(+)}$. 
 Moreover, the first author and R. Velasquez \cite{FggVr} 
proved the WPT for Jordan superalgebras of type $\josnm$. 
The WPT for Jordan superalgebras of types $\jpd$, Jordan superalgebra of Poisson bracket, 
$\mathcal{Q}_n$, and 
$\kap\oplus\cdots\oplus \kap\oplus\corf\cdot 1$ 
will be considered in other papers.

In the present paper the authors show that an analogue 
of the WPT holds when the simple quotient is of type $\jp$, $n\geq 3$. 
This paper has two sections. 
Section \ref{Preliminares} gives some preliminary results 
from the theory of Jordan superalgebras including those for the proof of the WPT 
and Section \ref{Main-Theorem} provides the proof of the WPT when the simple quotient is of type 
$\jp$, $n\geq 3$.

\section{Basic definitions and notation.}\label{Preliminares}

Throughout the paper, superalgebra means a finite-dimensional superalgebra over $\corf$.

\begin{definition}
An algebra $\alga$ is called a  superalgebra, if it is the direct sum of two non-zero subspaces 
$\alga_{\bar{0}}\oplus \alga_{\bar{1}}$ 
satisfying the multiplicative relation 
$\alga_{\bar{i}}\,\alga_{\bar{j}}\subseteq \alga_{\bar{i}+\bar{j} \,(\mathrm{mod}\, 2)}$. 
\end{definition}
 
An element $x\in\alga$ is called homogeneous of parity $\bar{i}$, 
if $x\in \alga_{\bar{i}}$, where $i=0,1$. We denote the parity of 
$x$, $x\in\alga_{\bar{i}}$, by 
$|x|=\bar{i}$, $i=0,1$.

\begin{definition}
A superalgebra $\jor$ is said to be a Jordan superalgebra,
if for all $x, y, z, t\in\jor_{\bar{0}} \dot\cup \jor_{\bar{1}}$ 
the superalgebra satisfies the superidentities
\begin{align}
&\label{idsupercom} xy=(-1)^{|x||y|}yx,\\
&\text{(supercommutative) and} \nonumber\\
&\label{idsupjordan}((xy)z)t+(-1)^{|t|(|z|+|y|)+|z||y|}((xt)z)y+(-1)^{|x|(|y|+|z|+|t|)+|t||z|}((yt)z)x\\
&\quad=(xy)(zt)+(-1)^{|t||z|+|t||y|}(xt)(yz)+(-1)^{|y||z|}(xz)(yt).\nonumber
\end{align}
\end{definition}

It is known that if 
$\alga$ is an associative superalgebra with multiplication 
$ab$, then $\alga^{(+)}$ is a Jordan superalgebra, 
where $\alga^{(+)}$ is a copy of the vector space $\alga$ 
and the multiplication in 
$\alga^{(+)}$ is defined by the 
supersymmetric product 
$a\circ b=\frac{1}{2}(ab+(-1)^{|a||b|}ba)$.
A Jordan superalgebra $\jor$ is called
\emph{special Jordan superalgebra}
if there exists an associative superalgebra 
$\alga$ such that 
$\jor\hookrightarrow\alga^{(+)}$.

For convenience, we denote by 
$ab$ the  product in associative superalgebras, 
$a\circ b$ the product in a special Jordan superalgebra, 
and 
$a\cdot b$ the product in arbitrary Jordan superalgebra.

Given an associative superalgebra 
$\alga$,  a 
graded linear mapping 
$\ast:\alga\longrightarrow\alga$ 
is said to be a \emph{superinvolution} in $\alga$, if 
$(a^\ast)^\ast=a$ and 
$(ab)^{\ast}=(-1)^{|a||b|}b^{\ast}a^{\ast}$ for every 
$a,b\in \alga_{\bar{0}}\,\dot\cup\,\alga_{\bar{1}}$.
If $\ast$ is a superinvolution in 
$\alga$, it is easy to see that 
the set 
$\displaystyle{\heraa=\{a\in\alga\,\mid\, a^\ast=a\}}$ 
is a subalgebra of $\alga^{(+)}$.

\medskip 
Note that 
$$\displaystyle{\begin{bmatrix}a& b \\ c& d \end{bmatrix}^\trp=
\begin{bmatrix} d^t& -b^ t \\ c^t& a^t\end{bmatrix} },$$ 
is a superinvolution in the associative superalgebra $\alMnn$. 
Then
$\mathcal{H}(\alMnn,\trp)$
is a Jordan superalgebra, 
that it is denoted  
 $\jp$. 
V. G. Kac \cite{Kac} proved that, i
f $n\geq 2$ then $\jp$ is a simple Jordan superalgebra.

Let $\jor$ be a finite-dimensional Jordan superalgebra and  let
$\bim=\bim_{\bar{0}}\oplus\bim_{\bar{1}}$  
be a 
$\jor$-bimodule. 
$\bim$ is called a Jordan bimodule, if the corresponding split null extension
$\mathcal{E}=\jor\oplus\bim$
is a Jordan superalgebra.
Recall that a split null extension is 
the algebra with underlying vector space equals to 
$\jor\oplus\bim$
with a multiplication that extends the multiplication in
$\jor$
through the action of $\jor$ on $\bim$,
while the product of two arbitrary elements in
$\bim$
is zero.

Let
$\bim$
be a
$\jor$-bimodule.
The opposite bimodule
$\obim=\obim_{\bar{0}}\oplus\obim_{\bar{1}}$
is defined by 
$\obim_{\bar{0}}=\bim_{\bar{1}},\,\obim_{\bar{1}}=\bim_{\bar{0}}$,
and  the action of
$\jor$ over $\obim$ by
$a\cdot m^{\textnormal{op}}=(-1)^{|a|}(am)^{\textnormal{op}}, \, m^{\textnormal{op}}\cdot a=(ma)^{\textnormal{op}}$
for all
$a\in\jor_{\bar{0}}\dot\cup\jor_{\bar{1}},\,m\in\obim_{\bar{0}}\dot\cup\obim_{\bar{1}}$.
Whenever
$\bim$
is a Jordan
$\jor$-superbimodule,
$\obim$ has a Jordan 
$\jor$-superbimodule structure.
Let $\bim=\jor$ as a vector superspace and let $a\cdot m$, $m\cdot a$ with $a\in\jor$, $m\in\bim$ 
be the products as defined in the superalgebra $\jor$.
 It is easy to see that $\bim$ has a natural structure of Jordan $\jor$-superbimodule. 
$\bim$ is called  the {\it{regular superbimodule}} of $\jor$ and it is denoted  by $\reg\,\jor$.

To study the general case of the WPT in Jordan superalgebra, 
the first author \cite{Fgg18} proved that
\begin{proposition}\label{pro2}
Let 
$\jor$ be a finite dimensional Jordan superalgebra with solvable radical
$\rad$ such that $\radd=0$ and 
$\jor/\rad$ is a semisimple Jordan superalgebra.
If the WPT holds for all simple quotient of Jordan superalgebra that appear in
$\jor/\rad$, then the WPT is valid for $\jor/\rad$.
\end{proposition} 

\begin{theorem}\label{teo: red2}
Let 
$\alga$
be a finite dimensional semisimple Jordan superalgebra, i.e 
$\rad(\alga)=0$,
where
$\rad$
is the solvable radical.
Let 
$\mathfrak{M}(\jor)$ be a class 
of finite dimensional Jordan
$\jor$-bimodules $\rad$ such that $\mathfrak{M}(\jor)$ 
is closed with respect to subbimodules and homomorphic images.
Denote by
$\mathcal{K}(\mathfrak{M},\jor)$
the class of finite dimensional Jordan superalgebras
$\alga$
that satisfy the following conditions:
\begin{enumerate}
\item
$\alga/\rad(\alga)\cong \jor$,
\item
$\radd(\alga)=0$,
\item
$\rad(\alga)\in\mj$.
\end{enumerate}
Then if the WPT is true for all superalgebras
$\algb\in\clasekmj$
with the restriction that the radical 
$\rad(\algb)$
is an irreducible
$\jor$-bimodule,
then it is true for all superalgebras
$\alga$
from
$\clasekmj$.
\end{theorem}

E. Zelmanov  \cite{Zelss} proved that 
$\jor/\rad$
is a semisimple Jordan superalgebra  if and only if
$\jor/\rad=(\jor_1\oplus\cdots\oplus\jor_k)\oplus(\kap\oplus\cdots\oplus\kap\oplus\corfd 1)$
 where $\jor_i$, $i=1,2\ldots,k$ is
 an unital simple Jordan superalgebra and
$\kap$ denotes the Kaplansky Jordan superalgebra. 
Using Zelmanov's theorem and Theorem \ref{teo: red2}, 
it is clear that  
verifying the analogue of the WPT for a Jordan superalgebra 
$\jor$ can be reduced to the proof of the WPT in the following cases
 $\jor/\rad$
 is simple unital;
 or 
 $\jor/\rad=(\mathcal{K}_3\oplus\mathcal{K}_3\oplus\cdots\oplus\kap)\oplus\corfd 1$.
Besides, simple Jordan superalgebras were classified by V. G. Kac \cite{Kac}, 
and I. L. Kantor \cite{Kan}. 
Observe that assuming $\radd=0$ and $\jor/\rad$ is a simple Jordan superalgebra,
 it is possible to reduce the problem  to consider only  Jordan irreducible  
$\jor/\rad$-superbimodules contained in $\rad$. 

Note that the irreducible Jordan superbimodules over the Jordan superalgebra 
$\jp$ were described  by C. Martinez and E. I. Zelmanov \cite{zelmar2}, 
who proved that only Jordan irreducible $\jp$-superbimodules
 are 
$\mathrm{Reg}\,(\jp)$, $\mathcal{P}_n$, and its opposites, where 
$\mathcal{P}_n=\{a\in\alMnn\,\mid\, a^{\mathrm{trp}}=-a\},$ $n\geq 3$.
Note that the action of $\jp$ over 
$\mathcal{P}_n$ 
is given by the supersymmetric product $a\circ m$ where $a\in\jp$, $m\in\mathcal{P}_n$.

By the results of C. Martinez, E. I. Zelmanov and I. P. Shestakov \cite{MarSheZel},
 we have that if $\jor$ is a Jordan superalgebra with unity 
$1$ and 
$\{e_1,\ldots,e_n\}\subseteq \jor_0$ 
is a set of pairwise orthogonal idempotents such that 
$1=\sum_{i=1}^n e_i$, 
then 
$\jor$ 
admits a Peirce decomposition
$$\jor=\Bigg(\bigoplus_{i=1}^n\jor_{ii}\Bigg)\bigoplus\Bigg(\bigoplus_{i<j}\jor_{ij}\Bigg), $$
where 
\begin{equation*}
\begin{aligned}
\jor_{ij}&=\left\{\,x\in\jor: \quad e_ix=\frac{1}{2}x,\quad e_jx=\frac{1}{2}x\,\right\}, \,\,\text{with}\,\, i\neq j, \\
\jor_{ii}&=\{\,x\in\jor: \quad e_ix=x\,\},
\end{aligned}
\end{equation*}
$\jor_{ii}$ and $\jor_{ij}$ 
are called the Peirce components of 
$\jor$ 
relative to the idempotents 
$e_i$, and $e_i, \,e_j$, respectively.
It is easy to see that the following relations hold
\begin{equation}\label{peircepro}
\begin{aligned}
&\jor_{ij}^2\subseteq \jor_{ii}+\jor_{jj},\quad \jor_{ij}\,\jor_{jk}\subseteq \jor_{ik}\text{ where } i\neq j, \, j\neq k,\\ 
&\jor_{ij}\,\jor_{kl}=0\quad \textnormal{where}\quad i\neq k,l; \quad j\neq k,l.
\end{aligned} 
\end{equation}

In the following part, we denote by 
$e_{ij}$
 the usual unit matrices. 
 For abbreviation, we write 
$\enn_{ij}=e_{n+i\,n+j}$, 
$\eun_{ij}=e_{i\,n+j}$, 
$\enu_{ij}=e_{n+i\,j}$ 
for $i,j=1,\ldots,n$.
We assume that 
\begin{equation*}
\begin{aligned}
u_{ij}=e_{ij}+\enn_{ji},\quad 
u_i=e_{ii}+\enn_{ii},\quad
s_{ij}=\eun_{ij}-\eun_{ji},\quad  
h_{ij}=\enu_{ij}+\enu_{ji},\quad
h_{i}=\enu_{ii}, \\
a_{i}=e_{ii}-e_{ii}^n,\quad  
a_{ij}=e_{ij}-e_{ji}^n,\quad  
b_{i}=e_{ii}^{1n}, \quad 
 b_{ij}=e_{ij}^{1n}+e_{ji}^{1n},\quad 
c_{ij}=e_{ij}^{n1}-e_{ji}^{n1}. 
\end{aligned}
\end{equation*}
With the above notation, 
it is easy to check that
$\{u_i,\,u_{ij},\, h_i,\,h_{ij} \}$
 is an additive basis of $\jp$ 
and 
$\{a_i,\,a_{ij},\, b_i,\,b_{ij}\}$
is an additive basis of $\mathcal{P}_n$
 for $i,j=1,\ldots, n$.


Note that 
$\{u_{1}, \ldots, u_{n}\}\subseteq(\jp)_0$ 
is a set of pairwise orthogonal idempotents in $\jp$, 
such that 
$u_1+\cdots+u_n=1$.
Thus, 
$\jp$ admits a Peirce decomposition with respect to 
$\{u_1,\ldots,u_n\}$.
Moreover, $(\jp)_{ii}$ is spanned by 
$\{u_i,h_i\}$,
while 
$(\jp)_{ij}$ is spanned by
$\{u_{ij},h_{ij},$ $s_{ij}\}$,
 and 
 \eqref{peircepro} holds for   
 $i,j=1,$ $\ldots,n$.

\section{Main theorem}\label{Main-Theorem}
In this section, we prove the main theorem of the present paper: 
the WPT holds when the simple quotient is of type $\jp$, $n\geq 3$.

\begin{theorem}\label{wptjpn}
Let $\jor$ be a finite dimensional Jordan superalgebra over 
a field 
$\corf$ of characteristic zero. 
Let $\rad$ be the solvable radical of $\jor$ such that $\radd=0$ and $\jor/\rad\cong\jp$, $n\geq 3$. 
Then, there exists a subsuperalgebra $\algs\subseteq\jor$ such that 
$\algs\cong\jp$ and $\jor=\algs\oplus\rad$.
\end{theorem}

Since $\jor$ is not necessarily a special Jordan superalgebra
we denote the product in $\jor$ by $a\cdot b$. 
The proof starts with the observation that there exists 
an additive basis
$\mathcal{B}=\{\wt{1},\wt{u}_i, \wt{u}_{ij}\mid i,j=1\ldots, n\}$
of $\algs_{\bar{0}}$, 
such that 
$\wt{a}\cdot \wt{b}= \widetilde{a\circ b}$ 
for all 
$\wt{a},\wt{b}\in\mathcal{B}$, $a,b\in(\jp)_{\bar{0}}$.
$\algs_{\bar{0}}\cong(\jp)_{\bar{0}}$ and 
$\jor_{\bar{0}}=\algs_{\bar{0}}\oplus\rad_{\bar{0}}$.

We need to prove that there exists 
an additive basis 
$\mathcal{B}^\prime=\{\wt{h}_{i},\wt{h}_{ij}, \wt{s}_{ij}\mid i,j=1\ldots, n\}$ 
of
 $\algs_{\bar{1}}\subseteq\jor_{\bar{1}}$ 
 such that 
 $\wt{a}\cdot \wt{b}=\widetilde{a\circ b}$
 for all $\wt{a},\wt{b}\in\mathcal{B}^\prime$, 
 $a,b\in(\jp)_{\bar{1}}$, 
 $\jor_{\bar{1}}=\algs_{\bar{1}}\oplus\rad_{\bar{1}}$.
Therefore, we obtain that there exists 
$\algs=\alga_{\bar{0}}\oplus\alga_{\bar{1}}$ 
such that 
$\algs\cong\jp$ and $\jor=\algs\oplus\rad$.

By Theorem \ref{teo: red2}, it suffices to consider $\rad$ isomorphic to
$\reg\jp$, $\mathcal{P}_n$, $(\reg\jp)^{\mathrm{op}}$ and  
$(\mathcal{P}_n)^{\mathrm{op}}$.
In each case, we can assume that 
$\sigma:\rad\longrightarrow \mathcal{M}$  
is the isomorphism between 
$\mathcal{N}$ and $\mathcal{M}$, 
where $\mathcal{M}$ is one of the irreducible $\jp$-superbimodules. 
Therefore,  we can assume that 
$\wt{a}\cdot {m}={\sigma^{-1}(a\circ \sigma(m)))}$, 
for all 
$\wt{a}\in\mathcal{B}\cup \mathcal{B}^\prime$,
${m}\in\rad$. 

Observe that
$\{\wt{u}_1, \ldots, \wt{u}_n\}$ 
 is a set of  pairwise orthogonal idempotents in $\jor_{\bar{0}}$, 
and  $\wt{1}=\wt{u}_1+\cdots+\wt{u}_n$, then  
$\jor$ admits Peirce decomposition.
Unless otherwise stated, we assume that 
$\wh_{i}$, $\wh_{ij}$, and $\ws_{ij}\in\jor_{\bar{1}}$ 
are  preimages of 
$\bh_i$, $\bh_{ij}$, and $\bs_{ij}$ 
under the canonical homomorphism
and 
$\overline{\mathcal{B}}=\{\bh_{i},\,\bh_{ij},\,\bs_{ij}\mid i,j=1\ldots, n \}$
is an additive basis of 
$\jor_{\bar{1}}/\rad_{\bar{1}}\cong(\jp)_{\bar{1}}$, 
such that 
$\bar{a}\cdot\bar{b}=\overline{a\circ b}$ 
for all $\bar{a}, \bar{b}\in\overline{\mathcal{B}}$, 
$a,b\in(\jp)_{\bar{1}}$.


Using Peirce properties, 
we conclude that 
$\wu_{i}\cdot\wh_i$, $\wu_{ji}\cdot\wh_{ji}\in (\jor_{\bar{1}})_{ii}$,
$\wu_{i}\cdot\wh_{ij}$, $\wu_{ij}\cdot\wh_{i}$, $\wu_{li}\cdot\wh_{lj}$, and $\wu_i\cdot\ws_{ij}\in (\jor_{\bar{1}})_{ij}$,
 while that 
 $\wh_{i}\cdot\ws_{ij}$, $\wh_{il}\cdot\ws_{lj}\in(\jor_{\bar{0}})_{ij}$ and $\wh_{ij}\cdot\ws_{ij}\in(\jor_{\bar{0}})_{i}+(\jor_{\bar{0}})_{j}$ for $i,j, l=1,\ldots, n$.
 
By \cite{zelmar2}, in order to prove our theorem, we have to consider 4 cases.

\medskip
{\bf Case 1. $\rad \cong\reg\jp$.} 
Assume that 
$\rad=\vect\langle e, v_i, v_{ij}\rangle\oplus \vect\langle g_i,g_{ij},z_{ij}\rangle $, 
where 
$1\leftrightarrow e$,
$u_i\leftrightarrow v_i$, 
$u_{ij}\leftrightarrow v_{ij}$,
$h_i\leftrightarrow g_i$,
$h_{ij}\leftrightarrow g_{ij}$,
$s_{ij}\leftrightarrow z_{ij}$ for $i,j=1,\ldots, n$.

Note that 
$(\rad_{\bar{0}})_{ij}$, $(\rad_{\bar{0}})_{i}$, $(\rad_{\bar{1}})_{i}$, and 
$(\rad_{\bar{1}})_{ij}$ are 
spanned  by 
$\{v_{ij},v_{ji}\}$, $\{v_{i}\}$, $\{g_{i}\}$, and $\{g_{ij},z_{ij}\}$, for $i,j=1,\ldots, n$, 
respectively. T
hen, we assume that there exist scalars 
$\eta_{i}$, $\eta_{ij}$, $\eta_{iij}$, $\eta_{iji}$, $\eta_{ijil}$, $\alpha_{iij}$, $\alpha_{iji}$, $\alpha_{ij}$, 
$\beta_{iij}$, $\beta_{ijjl}$, $\gamma_{iij}$, $\gamma_{ijjl}$, 
$\Lambda^{i}_{ij}$, $\Lambda^{j}_{ij}$, $\Lambda^{ji}_{iij}$ $\Lambda^{ij}_{iij}$, 
$\Lambda_{ijjl}^{il}$, and $\Lambda_{ijjl}^{li}$ such that 
$\wu_i \cdot\wh_i=\wh_i+\eta_{i} g_i,$ $u_i \cdot \wh_{ij}=\frac{1}{2}\wh_{ij}+\eta_{iij}g_{ij}+\alpha_{iij}z_{ij},$
$\wu_{ji} \cdot \wh_{ij}=\wh_i+\eta_{ji}g_i$, 
$\wu_{ij}\cdot  \wh_i=\frac{1}{2}\wh_{ij}+\eta_{iji}g_{ij}+\alpha_{iji}z_{ij}$,  
$\wu_{ij} \cdot \wh_{il}=\frac{1}{2}\wh_{jl}+\eta_{ijil}g_{jl}+\alpha_{ijil}z_{jl}$, 
$\wu_{i}\cdot  \ws_{ij}=\frac{1}{2}\ws_{ij}+\gamma_{iij}g_{ij}+\beta_{iij}z_{ij}$,  
$\wu_{ij}\cdot  \ws_{jl}=\frac{1}{2}\ws_{il}+\gamma_{ijjl}g_{il}+\beta_{ijjl}z_{il}$, 
$\wh_{i}\cdot \ws_{ij}=\frac{1}{2}\wu_{ji}+\Lambda_{iij}^{ji}v_{ji}+\Lambda^{ij}_{iij}v_{ij}$,
$\wh_{ij}\cdot\ws_{jl}=\frac{1}{2}\wu_{li}+\Lambda_{ijjl}^{il}v_{il}+\Lambda_{ijjl}^{li}v_{li}$,  
and 
$\wh_{ij}\cdot\ws_{ij}=\frac{1}{2}(\wu_{j}-\wu_{i})+\Lambda_{ij}^{i}v_{i}+\Lambda_{ij}^{j}v_{j}$.

\begin{lemma}\label{lemauihihijsij}
$\displaystyle{\wu_i\cdot\wh_i=\wh_i, \quad \wu_i \cdot\wh_{ij}=\frac{1}{2}\wh_{ij} \quad\text{and}\quad
\wu_{i}\cdot \ws_{ij}= \frac{1}{2}\ws_{ij}.}$
\end{lemma}
\proof  Using \eqref{idsupjordan}, 
we proceed to determine constants 
$\eta_{i}$, $\eta_{ij}$, $\eta_{iij}$, $\eta_{iji}$, $\eta_{ijil}$, $\alpha_{iij}$, $\alpha_{iji}$, $\alpha_{ij}$, 
$\beta_{iij}$, $\beta_{ijjl}$, $\gamma_{iij}$, $\gamma_{ijjl}$, 
$\Lambda^{i}_{ij}$, $\Lambda^{j}_{ij}$, $\Lambda^{ji}_{iij}$ $\Lambda^{ij}_{iij}$, 
$\Lambda_{ijjl}^{il}$, and 
$\Lambda_{ijjl}^{li}\in\corf$.

Now, replacing 
$x=z=t$ by $\wu_i$ and $y$ by $\wh_i$ in \eqref{idsupjordan},
we obtain that
\begin{equation}\label{id: uihi}
2((\wu_i\cdot \wh_i)\cdot \wu_i)\cdot \wu_i+\wu_i\cdot\wh_i=3\wu_i\cdot(\wu_i\cdot\wh_i).
\end{equation}
By the multiplication and the action of 
$\mathcal{J}$ over $\mathcal{N}$ in \eqref{id: uihi}, 
we get $7\eta_{i} g_i=6\eta_{i} g_i$. Thus, 
$\eta_{i}=0$, and $\wu_i \cdot \wh_i=\wh_i$. 
In a similar way, changing $\wh_{i}$ by $\wh_{ij}$ in \eqref{id: uihi},
implies that 
\begin{equation}\label{id: uihij}
2((\wu_i\cdot \wh_{ij})\cdot \wu_i)\cdot \wu_i+\wu_i\cdot\wh_{ij}=3\wu_i\cdot(\wu_i\cdot\wh_{ij}).
\end{equation}
Thus,
$\displaystyle{\eta_{iij}g_{ij}+\alpha_{iij}z_{ij}=0}$. This equality implies that 
$\eta_{iij}=\alpha_{iij}=0$. Therefore, we conclude that 
$\wu_i \cdot \wh_{ij}=\frac{1}{2}\wh_{ij}$.
Analogously, we obtain 
$\gamma_{iij}=\beta_{iij}=0$. We get that
$\wu_i \cdot \ws_{ij}=\frac{1}{2}\ws_{ij}$.
\endproof
\begin{lemma}\label{lemappal} 
There exist $\xi_{ij}\in\corf$, 
$i,j=1,\ldots,n$ such that
\begin{align}
\begin{aligned}\label{J01reg}
\wu_{ji} \cdot \wh_{ij}&=\wh_i, +\xi_{ji}g_i, &&   
\wu_{ij}\cdot  \wh_i=\frac{1}{2}\wh_{ij}+\frac{1}{2}\xi_{ji}g_{ij}, \\ 
\wu_{ij} \cdot \wh_{il}&=\frac{1}{2}\wh_{jl}+\frac{1}{2}(\xi_{il}-\xi_{jl})g_{jl}, &&
\wu_{ij}\cdot  \ws_{jl}=\frac{1}{2}\ws_{il}+\frac{1}{2}(\xi_{il}-\xi_{jl})z_{il}.
\end{aligned}\\
\begin{aligned}\label{J11reg}
\wh_{i}\cdot \ws_{ij}&=\frac{1}{2}\wu_{ji}+\frac{1}{2}\xi_{ji}v_{ji},\quad
\wh_{ij}\cdot\ws_{jl}=\frac{1}{2}\wu_{il}+\frac{1}{2}(\xi_{ij}-\xi_{lj})v_{il},\\
\wh_{ij}\cdot\ws_{ij}&=\frac{1}{2}(\wu_{j}-\wu_{i}).
\end{aligned}\
\end{align}
Where \eqref{J01reg} are commutative products,
and
\eqref{J11reg} are anti-commutative products.
\end{lemma}
\proof We start by showing that 
$\alpha_{iji}=\alpha_{ijil}=\gamma_{ijjl}=\Lambda_{iij}^{ij}=\Lambda_{ijjl}^{il}=0$,
and $\Lambda_{ij}^i=-\Lambda_{ij}^j$.
Substituting 
$x=t$ by $\wu_{ij}$ and $y=z$ by $\wh_i$ 
($y=z$ by $\wh_{il}$ and $y=z$ by $\ws_{jl}$, respectively)
 in  \eqref{idsupjordan}, we obtain that
 $$0=((\wu_{ij}\cdot\wh_i)\cdot\wh_i)\wu_{ij}=\alpha_{iji}(z_{ij}\cdot\wh_i)\wu_{ij}=-\alpha_{iji}v_{ji}\cdot\wu_{ij}.$$
 Then, we conclude that 
 $\alpha_{iji}=0$ ($\alpha_{ijil}=0$ and $\gamma_{ijjl}=0$, respectively).
Replacing
 $x$ by $\wh_{i}$, $y$ by $\ws_{ij}$, $z$ by $\wu_{ij}$ and $t$ by $\wu_{ji}$ in \eqref{idsupjordan}, 
 we get
 \begin{equation}\label{id: hisijreg}
((\wh_i\cdot\ws_{ij})\cdot\wu_{ij})\cdot\wu_{ji}=\frac{1}{2}(\wu_i+\wu_j)\cdot(\wh_i\cdot\ws_{ij}).
 \end{equation}
From \eqref{id: hisijreg}, it is clear that 
 $\frac{1}{2}\Lambda_{iij}^{ji}v_{ji}=\frac{1}{2}\Lambda_{iij}^{ij}v_{ij}+\frac{1}{2}\Lambda_{iij}^{ji}v_{ji}$. Then, we conclude that $\Lambda_{iij}^{ij}=0$. 
 We denote $\Lambda_{iij}^{ji}=\Lambda_{iij}$, 
 thus we write
 $\displaystyle{\wh_i\cdot\ws_{ij}=\displaystyle{\frac{1}{2}\wu_{ji}+\Lambda_{iij}v_{ji}}}$.
Similarly, we obtain $\Lambda_{ijjl}^{il}=0$ and  
 $\displaystyle{\wh_{ij}\cdot\ws_{jl}=\displaystyle{\frac{1}{2}\wu_{li}+\Lambda_{ijjl}v_{li}}}$,
 where 
$\Lambda_{ijjl}^{li}=\Lambda_{ijjl}$.

Now, substituting 
 $x$ by $\wh_{ij}$, $y$ by $\ws_{ij}$, $z$ by $\wu_{ij}$ and $t=\wu_{ji}$ in \eqref{idsupjordan}, 
 we obtain 
 \begin{equation}\label{id: hijsij}
((\wh_{ij}\cdot\ws_{ij})\cdot\wu_{ij})\cdot\wu_{ji}+((\wu_{ji}\cdot\wh_{ji})\cdot\wu_{ij})\cdot\ws_{ij}=
 (\wh_{ij}\cdot\ws_{ij})\cdot(\wu_{ij}\cdot\wu_{ji}).
 \end{equation}
It is easy to see that calculations on the left and right of \eqref{id: hijsij} give
$$(\Lambda_{ij}^i+\Lambda_{ij}^i-(2\eta_{iji}+\eta_{ji}))v_i+(\Lambda_{ij}^i+\Lambda_{ij}^i+(2\eta_{iji}+\eta_{ji}))v_j=0.$$
By the linear independence of  $v_i$ and $v_j$, we conclude that
 $\Lambda_{ij}^i+\Lambda_{ij}^j=0$ and 
\begin{equation}\label{ec:ejieiji}
\eta_{ji}+2\eta_{iji}=0.
\end{equation}
We denote $\Lambda_{ij}^j=\Lambda_{ij}$. Thus, we write
  $\displaystyle{ \wh_{ij}\cdot\ws_{ij}=\frac{1}{2}(\wu_{j}-\wu_i)+\Lambda_{ij}(v_j-v_i)}$.

Further, substituting
$x$ by $t=\wu_{ij}$, $y$ by $\wh_{il}$ and $z$ by $\wu_{ji}$  in  \eqref{idsupjordan},  
we get
$((\wu_{ij}\cdot\wh_{il})\cdot\wu_{ji})\cdot\wu_{ij}=(\wu_{ij}\cdot\wh_{il})\cdot(\wu_{ji}\cdot\wu_{ij})$
 and simplifying this equality, we obtain 
\begin{equation}\label{ec:eijlejil}
\eta_{ijil}+\eta_{jijl}=0.
\end{equation}

Also, replacing
$x$ by $\wu_{i}$, $y$ by $\wh_{ij}$, $z$ by $\wu_{il}$ and 
$t$ by $\wu_{lj}$ in  \eqref{idsupjordan}, 
we  get
$((\wu_{i}\cdot\wh_{ij})\cdot\wu_{il})\cdot\wu_{lj}=(\wu_{i}\cdot\wh_{ij})\cdot(\wu_{il}\cdot\wu_{lj})$.
By the linear independence, we conclude that $\eta_{ij}-\eta_{lj}=2\eta_{ilij}$. 
Combining this equality with \eqref{ec:ejieiji}, we get
\begin{equation}\label{ec:eijeljeilij}
\frac{1}{2}(\eta_{ij}-\eta_{lj})=\eta_{ilij}=\eta_{jlj}-\eta_{jij}.
\end{equation}

Now, substituting 
 $x=t$ by $\wu_{ij}$, $y$ by $\ws_{jl}$ and $z$ by $\wu_{ji}$ in \eqref{idsupjordan}, 
we obtain 
$((\wu_{ij}\cdot\ws_{jl})\cdot\wu_{ji})\cdot\wu_{ij}=\frac{1}{2}(\wu_{ij}\cdot\ws_{jl})\cdot(\wu_j+\wu_i).$
From this equality, an easy computation gives that
\begin{equation}\label{ec:bijlbjil}
\beta_{ijjl}+\beta_{jiil}=0. 
\end{equation}
Similarly, we obtain that
$((\wu_{ij}\cdot\ws_{jl})\cdot\wu_{jl})\cdot\wu_{li}+\frac{1}{4}\ws_{jl}=
\frac{1}{2}(\wu_{ij}\cdot\ws_{jl})\cdot\wu_{ji},$
combining this with \eqref{ec:bijlbjil},
 it follows that
\begin{equation}\label{ec:blijbjlibijl}
\beta_{ijjl}+\beta_{jlli}+\beta_{liij}=0.
\end{equation}

Changing 
$x$ by $\wu_{i}$, $y$ by $\wh_{i}$, $z$ by $\wu_{ij}$ and $t$ by $\ws_{ij}$ in \eqref{idsupjordan}, 
we obtain
\begin{equation}\label{id: hisij}
(\wu_{ij}\cdot\wh_{i})\cdot\ws_{ij}+((\wh_i\cdot\ws_{ij})\cdot\wu_{ij})\cdot\wu_{i}=
\frac{1}{2}((\wu_{ij}\cdot\wh_i)\cdot\ws_{ij}+\wu_{ij}\cdot(\wh_i\cdot\ws_{ij})).
\end{equation} 
By the linear independence, we get
\begin{equation}\label{ec:LiijLijeijj}
\Lambda_{iij}=\Lambda_{ij}+\eta_{iji}.
\end{equation}

Now, replacing
$x$ by $\wu_{ij}$, $y$ by $\ws_{jl}$, $z$ by 
$\wh_{i}$ and $t$ by $\ws_{ij}$ in \eqref{idsupjordan}, 
we assert that
\begin{equation}\label{ec: uijsjlhisij}
((\wu_{ij}\cdot\ws_{jl})\cdot\wh_{i})\cdot\ws_{ij}=(\wu_{ij}\cdot\ws_{jl})\cdot(\wh_{i}\cdot\ws_{ij}).
\end{equation} 
And simplifying  \eqref{ec: uijsjlhisij}, it follows that
$\beta_{liij}-\beta_{jiil}=\Lambda_{iij}-\Lambda_{iil}.$ 
Besides, combining \eqref{ec:bijlbjil}, \eqref{ec:blijbjlibijl} and this equality, 
it is easy to see that 
 \begin{equation}\label{ec:bljiLiijLiil}
 \beta_{jlli}=\Lambda_{iil}-\Lambda_{iij}.
 \end{equation}
 
 Now, substituting 
 $x$ by $\wu_{ij}$, $y$ by $\ws_{jl}$, $z$ by $\wh_{i}$, and
  $t$ by $\wh_{l}$ in \eqref{idsupjordan}, 
 we get
\begin{equation}\label{id: sjlhihl}
((\wu_{ij}\cdot\ws_{jl})\cdot\wh_i)\cdot \wh_l=-(\wh_{l}\cdot\ws_{lj})\cdot(\wu_{ij}\cdot\wh_i).
\end{equation}
Using \eqref{id: sjlhihl}, it is easy to see that 
 \begin{equation}\label{ec:ELBA}
 \eta_{lil}+\Lambda_{iil}+\beta_{ijjl}=\eta_{jlji}+\eta_{iji}+\Lambda_{llj}.
 \end{equation}
 
Substituting \eqref{ec:bljiLiijLiil} into \eqref{ec:ELBA},
 we conclude that
$\Lambda_{iil}-\Lambda_{lli}=\eta_{iji}-\eta_{lil}+\eta_{jlji}$
holds. Again, replacing \eqref{ec:eijeljeilij} in the above equality, we obtain 
$\Lambda_{iil}-\Lambda_{lli}=\eta_{ili}-\eta_{lil}$. 
Now, replacing \eqref{ec:LiijLijeijj} in the above equality, it is clear that 
$\Lambda_{il}=\Lambda_{li}.$
Besides, observe that 
$\wh_{li}\cdot\ws_{li}=-\wh_{il}\cdot\ws_{il}$, then
this gives that
$\Lambda_{il}=-\Lambda_{li}$. So, we conclude that
$\Lambda_{il}=0$. Finally, we write 
$\wh_{ij}\cdot\ws_{ij}=\frac{1}{2}(\wu_j-\wu_i).$

Further, writing $\Lambda_{il}=0$ in \eqref{ec:LiijLijeijj}, 
we conclude  that
$\Lambda_{iij}=\eta_{iji}$. 
Now, replacing this equality in \eqref{ec:bljiLiijLiil} 
and using \eqref{ec:eijeljeilij}, 
we deduce that
\begin{equation}\label{ec:e=b}
\beta_{jlli}=\eta_{jlji}.
\end{equation}

Substituting  
 $x$ by $\wu_{l}$, $y$ by $\wh_{l}$, $z$ by $\ws_{lj}$ and $t$ by $\wh_{ji}$ in \eqref{idsupjordan}, 
 we obtain
\begin{equation}\label{id: hlsljhij}
((\wu_l\cdot\wh_l)\cdot\ws_{lj})\cdot\wh_{ij}=(\wu_l\cdot\wh_l)\cdot(\ws_{lj}\cdot\wh_{ij}).
\end{equation}
From \eqref{id: hlsljhij} and the linear independence, we give
$\Lambda_{ijjl}=\eta_{jlji}+\Lambda_{llj}-\eta_{lil}=\eta_{jlji}+\eta_{ljl}-\eta_{lil}=\eta_{jlji}+\eta_{ijil}.$
Now, using 
\eqref{ec:blijbjlibijl}  and \eqref{ec:e=b}
in these equalities, we obtain
\begin{equation}\label{ecuacion}
\Lambda_{jiil}=\beta_{jlli}=\eta_{jlji}=\frac{1}{2}(\eta_{ji}-\eta_{li}).
\end{equation}

Finally, with \eqref{ecuacion}, we have proved that:
\begin{equation}\label{ecuacion1}
\begin{aligned}
\wu_{ij}\cdot\wh_{ij}&=\wh_j+\eta_{ij}g_j, &&
\wu_{ij}\cdot\wh_{i}=\frac{1}{2}\wh_{ij}-\frac{1}{2}\eta_{ij}g_{ij},  \\
\wu_{ij}\cdot\wh_{il}&=\frac{1}{2}\wh_{jl}+\frac{1}{2}(\eta_{il}-\eta_{jl})g_{jl}, &&
\wu_{ij}\cdot\ws_{jl}=\frac{1}{2}\ws_{il}+\frac{1}{2}(\eta_{il}-\eta_{jl})z_{il}, \\
\wh_{i}\cdot\ws_{ij}&=\frac{1}{2}\wu_{ji}-\eta_{ji}v_{ji}, &&
\wh_{ij}\cdot\ws_{jl}=\frac{1}{2}\wu_{li}+\frac{1}{2}(\eta_{ij}-\eta_{lj})v_{li}.
\end{aligned}
\end{equation}
Making $\xi_{ij}=\eta_{ij}$ in \eqref{ecuacion1}, the proof of  the lemma is completed.
\endproof

\medskip\noindent 
{\bf Case 2. $\rad\cong(\reg\jp)^{\textrm{op}}$.}
\begin{lemma}\label{le: regop01}
$\wu_i\cdot\wh_{i}= \wh_i$,  
$\wu_{i}\cdot \wh_{ij}=\frac{1}{2}\wh_{ij}, $
$\wu_{i}\cdot \ws_{ij}=\frac{1}{2}\ws_{ij}, $
$\wu_{ij}\cdot\wh_{il}=\frac{1}{2}\wh_{jl}, $
$\wu_{ij}\cdot\ws_{jl}=\frac{1}{2}\ws_{il}, $
$\wu_{ij}\cdot\wh_{ij}=\wh_{j},$  
$\wu_{ij}\cdot\wh_{i}=\frac{1}{2}\wh_{ij}$,
$\wh_{ij}\cdot\ws_{ij}=\frac{1}{2}(\wu_j-\wu_i)$,
$\wh_{i}\cdot\ws_{ij}=\frac{1}{2}\wu_{ji}$
and  
$\wh_{ij}\cdot\ws_{jl}=\frac{1}{2}\wu_{li}.$

\end{lemma}
\proof 
Analogously to the proof of Lemma \ref{lemauihihijsij},   
we prove the equalities 
$\wu_i\cdot\wh_i=\wh_i$,
$\wu_i\cdot\wh_{ij}=\frac{1}{2}\wh_{ij}$ and
$\wu_i\cdot\ws_{ij}=\frac{1}{2}\ws_{ij}$.
 
Now, assume that 
$\wu_{ij}\cdot \wh_i=\frac{1}{2}\wh_{ji}+\eta_{iji}^{ji}v_{ji}+\eta_{iji}^{ij}v_{ij}$ 
and $\wu_{ij}\cdot \wh_{ij}=\wh_j+\eta_{ij}v_j$ for some 
$\eta_{iji}^{ji}$, $\eta_{iji}^{ij}$ and $\eta_{ij}\in\corf$.
Substituting $x=t$ by $\wu_{ij}$ and $y=z$ by $\wh_i$ in \eqref{idsupjordan},
we get $((\wu_{ij}\cdot \wh_i)\cdot\wh_i)\cdot\wu_{ij}=0$.
Then, it follows that 
$\eta_{iji}^{ij}=0$. If we replace
$x=z=t$ by $\wu_{ij}$ and $y$ by $\wh_i$ in \eqref{idsupjordan}, we obtain that
$((\wu_{ij}\cdot\wh_i)\cdot \wu_{ij})\cdot \wu_{ij}=0$.
Therefore, we conclude that 
$\frac{1}{2}\eta_{ij}+\eta_{iji}^{ji}=0$.
Changing 
$x$ by $\wu_{ij}$, $y$ by $\wh_i$, $z$ by $\wh_j$ and $t$ by 
$\wu_{ji}$ in \eqref{idsupjordan},  
we get
\begin{equation}\label{id: uijhiluiluli}
((\wu_{ij}\cdot\wh_i)\cdot\wh_j)\cdot\wu_{ji}=(\wu_{ij}\cdot\wh_i)\cdot (\wu_{ij}\cdot\wh_i).
\end{equation}
Simplifying  \eqref{id: uijhiluiluli} and using the linear independence of the elements, 
we have that $\eta_{iji}^{ji}=0$. This implies that $\eta_{ij}=0$ for all $i,j=1,\ldots,n$.
So, we write
$\wu_{ij}\cdot \wh_i=\frac{1}{2}\wh_{ji}$ and 
$\wu_{ij}\cdot \wh_{ij}=\wh_{j}$. 
Further, substituting 
$x$ by $\wu_{ij}$, $y$ by $\wh_{il}$, $z$ by $\wu_{il}$ and $t$ by $\wu_{li}$ 
in \eqref{idsupjordan}, it follows easily that $\wu_{ij}\cdot\wh_{il}=\frac{1}{2}\wh_{jl}$.

Let
$\wu_{ij}\cdot\ws_{jl}=\frac{1}{2}\ws_{il}+\beta_{ijjl}^{il}v_{il}+\beta_{ijjl}^{li}v_{li}$, 
where $\beta_{ijjl}^{il}$ and $\beta_{ijjl}^{li}\in\corf$. 
In the same manner, observe that $((\wu_{ij}\cdot\ws_{jl})\cdot\ws_{ij})\cdot\wu_{ij}=0$. 
Thus, we obtain 
$\beta_{ijjl}^{li}=0$. Similarly, substituting 
$x$ by $\wu_{ij}$, $y=z$ by $\ws_{jl}$ and $t$ by $\wu_{jl}$ in \eqref{idsupjordan}, it is clear that 
$\beta_{ijjl}^{il}=0$. 
Then, we conclude that 
$\wu_{ij}\cdot\ws_{jl}=\frac{1}{2}\ws_{il}$. 

Now, assuming that 
$\wh_{ij}\cdot\ws_{ij}=\frac{1}{2}(\wu_j-\wu_i)+\Lambda_{ij}^ig_i+\Lambda_{ij}^jg_j$,
$\wh_{i}\cdot\ws_{ij}=\frac{1}{2}\wu_{ji}+\Lambda_{iij}^gg_{ji}+\Lambda_{iij}^zz_{ji}$, 
and 
$\wh_{ij}\cdot\ws_{jl}=\frac{1}{2}\wu_{li}+\Lambda_{ijjl}^gg_{li}+\Lambda_{ijjl}^zz_{li}$
and using \eqref{id: hijsij}, we get $\Lambda_{ij}^i=0$. 
Substituting 
$x$ by $\wh_{ij}$, $y$ by $\ws_{ij}$, $z$ by $\wu_{ji}$ and $t$ by $\wu_{ij}$ 
in \eqref{idsupjordan}, 
it is clear that $\Lambda_{ij}^j=0$. 
Besides, replacing 
$x$ by $\wh_{i}$, $y$ by $\ws_{ij}$, $z$ by $\wu_{ij}$ and $t$ by $\wu_{ji}$ in 
\eqref{idsupjordan}, 
we deduce that $\Lambda_{iij}^z=0$. Similarly, 
$\Lambda_{jji}^g=\Lambda_{ijjl}^z=\Lambda_{ijjl}^g=0$, which proves the lemma.
\endproof
 
{\bf Case 3. $\rad\cong\mathcal{P}_n$.} 
Let 
$\rad=\vect\langle w_i, w_{ij}\rangle\oplus \vect\langle y_i,y_{ij},x_{ij}\rangle$
and 
$\sigma({w}_i)=a_{i}$,  $\sigma(w_{ij})=a_{ij}$, 
$\sigma(y_i)=b_i$, 
$\sigma(y_{ij})=b_{ij}$, and 
$\sigma(x_{ij})=c_{ij}$. 
It is clear that 
 $(\rad_{\bar{0}})_{ij}$, $(\rad_{\bar{0}})_{i}$, $(\rad_{\bar{1}})_{i}$ and $(\rad_{\bar{1}})_{ij}$ 
 are spanned  by 
$\{w_{ij},w_{ji}\}$, $\{w_{i}\}$, $\{y_{i}\}$ and $\{y_{ij},x_{ij}\}$, respectively.

Now, Let 
$\eta_{i}$, $\eta_{ij}$, $\eta_{iij}$, $\eta_{ijil}$, $\eta_{iji}$, $\alpha_{iij}$, $\alpha_{iji}$, $\alpha_{ij}$, 
$\beta_{iij}$, $\beta_{ijjl}$, $\gamma_{iij}$, $\gamma_{ijjl}$, 
$\Lambda^{i}_{ij}$, $\Lambda^{j}_{ij}$, $\Lambda^{ji}_{iij}$, $\Lambda^{ij}_{iij}$, 
$\Lambda_{ijjl}^{il}$ and $\Lambda_{ijjl}^{li}$  scalars in $\corf$
such that 
$\wu_i\cdot\wh_i=\wh_i+\eta_{i} y_i,$ 
$\wu_i\cdot\wh_{ij}=\frac{1}{2}\wh_{ij}+\eta_{iij}y_{ij}+\alpha_{iij}x_{ij},$
$\wu_{ji}\cdot\wh_{ij}=\wh_i+\eta_{ji}y_i$, 
$\wu_{ij}\cdot\wh_i=\frac{1}{2}\wh_{ij}+\eta_{iji}y_{ij}+\alpha_{iji}x_{ij}$,  
$\wu_{ij}\cdot\wh_{il}=\frac{1}{2}\wh_{jl}+\eta_{ijil}y_{jl}+\alpha_{ijil}x_{jl}$, 
$\wu_{i}\cdot\ws_{ij}=\frac{1}{2}\ws_{ij}+\gamma_{iij}y_{ij}+\beta_{iij}x_{ij}$,  
$\wu_{ij}\cdot\ws_{jl}=\frac{1}{2}\ws_{il}+\gamma_{ijjl}y_{il}+\beta_{ijjl}x_{il}$, 
$\wh_{i}\cdot\ws_{ij}=\frac{1}{2}\wu_{ji}+\Lambda_{iij}^{ji}w_{ji}+\Lambda^{ij}_{iij}w_{ij}$,
$\wh_{ij}\cdot\ws_{jl}=\frac{1}{2}\wu_{li}+\Lambda_{ijjl}^{il}w_{il}+\Lambda_{ijjl}^{li}w_{li}$,  
$\wh_{ij}\cdot\ws_{ij}=\frac{1}{2}(\wu_{j}-\wu_{i})+\Lambda_{ij}^{i}w_{i}+\Lambda_{ij}^{j}w_{j}$, for $i,j, l=1,\ldots,n$.

\begin{lemma}\label{le: skew01}
$\wu_i\cdot\wh_{i}= \wh_i$, 
$\wu_{i}\cdot\wh_{ij}=\frac{1}{2}\wh_{ij}$, 
$\wu_{i}\cdot\ws_{ij}=\frac{1}{2}\ws_{ij}$, 
$\wu_{ij}\cdot\wh_{il}=\frac{1}{2}\wh_{jl}$, 
$\wu_{ij}\cdot\ws_{jl}=\frac{1}{2}\ws_{il}$, 
$\wu_{ij}\cdot\wh_{ij}=\wh_{j}$, 
$\wu_{ij}\cdot\wh_{i}=\frac{1}{2}\wh_{ij}$, 
$\wh_{ij}\cdot\ws_{ij}=\frac{1}{2}(\wu_j-\wu_i)$,  
$\wh_j\cdot\ws_{jl}=\frac{1}{2}\wu_{lj}$,
$\wh_{ji}\cdot\ws_{il}=\frac{1}{2}\wu_{lj}$.
\end{lemma}
\proof
We can now proceed analogously to the proof of Lemma \ref{lemauihihijsij}, 
and obtain that 
$\wu_i\cdot\wh_{i}= \wh_i,$
$\wu_{i}\cdot \wh_{ij}=\frac{1}{2}\wh_{ij}$ and finally that
$\wu_{i}\cdot \ws_{ij}=\frac{1}{2}\ws_{ij}$.
Similarly to the proof of Lemma \ref{lemappal}, we conclude that 
$\eta_{iji}=\eta_{ijil}=\beta_{ijjl}=0$. 

Now, we proceed to finding the constants. 
Substitute 
$x=t$ by $\wu_{ij}$, $y$ by $\wh_{il}$ and $z$ by $\wu_{ji}$ 
($x=t$ by $\wu_{ij}$, $y$ by $\ws_{jl}$ and $z$ by $\wu_{ji}$, respectively) 
in \eqref{idsupjordan}. 
It follows easily that 
$\alpha_{ijil}=0$, 
($\gamma_{ijjl}=0$, respectively). 
So, we write 
$\displaystyle{\wu_{ij}\cdot\wh_{il}=\frac{1}{2}\wh_{jl}}$ and 
$\displaystyle{\wu_{ij}\cdot\ws_{jl}=\frac{1}{2}\ws_{il}}$. 
Moreover, replacing 
$x=t$ by $\wu_{ij}$, $y$ by $\wh_{ij}$ and $z$ by $\wu_{ji}$ in 
\eqref{idsupjordan}, 
we deduce that $\eta_{ij}=0$. 
Thus, $\wu_{ij}\cdot\wh_{ij}=\wh_j$.

Changing $x=t$ by $\wu_{ij}$, $y$ by $\wh_{i}$ and $z$ by $\wu_{jl}$ in 
\eqref{idsupjordan}, we obtain
$\alpha_{iji}x_{jl}=\alpha_{iji}x_{lj}$. 
Observe that $x_{lj}=-x_{jl}$ and by the linear independence of $x_{lj}$, 
we conclude that 
$\alpha_{iji}=0$.
 Finally, we write 
 $\wu_{ij}\cdot\wh_i=\frac{1}{2}\wh_{ij}$.

By \eqref{id: hisijreg}, it follows that $\Lambda_{iij}^{ij}=0$.
Similarly, we obtain  $\Lambda_{ijjl}^{il}=0$.
Now, by \eqref{id: hijsij} we give 
$\Lambda_{ij}^i+\Lambda_{ij}^j=0$. 
Denote $\Lambda_{ij}=\Lambda_{ij}^j$; 
then we write 
$\wh_{ij}\cdot\ws_{ij}=\frac{1}{2}(\wu_j-\wu_i)+\Lambda_{ij}(w_j-w_i)$.
Further, by \eqref{id: hisij}, we obtain  
$\Lambda_{ij}=\Lambda_{iij}$. 
By \eqref{ec: uijsjlhisij} it
follows that 
$\Lambda_{iij}=-\Lambda_{iil}$.
Now, by \eqref{id: sjlhihl} it is clear that 
$\Lambda_{llj}=-\Lambda_{iil}$,  
thus we conclude that
$\Lambda_{iij}=\Lambda_{llj}$, 
which is the same as
$\Lambda_{lj}=\Lambda_{ij}$.
Replace
$x$ by $\wu_{il}$, $y$ by $\ws_{lj}$, $z$ by $\wu_{li}$ and $t$ by $\wh_{jl}$ 
in \eqref{idsupjordan}. 
An easy calculation shows that 
$\Lambda_{jl}=\Lambda_{ij}$. 
Combining these last equations, we conclude that 
$\Lambda_{jl}=\Lambda_{ij}=\Lambda_{lj}$. 
Furthermore, note that 
$\Lambda_{lj}=-\Lambda_{jl}=\Lambda_{jl}$, 
then
$\Lambda_{llj}=\Lambda_{lj}=0$. Finally,
by \eqref{id: hlsljhij} it is clear that $\Lambda_{ijjl}=\Lambda_{llj}=0$,
for all $i,\,j,\,l=1,\ldots,n$, and the proof is completed.
\endproof

{\bf{Case 4.} $\rad\cong(\mathcal{P}_n)^{\mathrm{op}}$.}
\begin{lemma}\label{le: skewop01}
$\wu_i\cdot\wh_{i}= \wh_i$,  
$\wu_{i}\cdot\wh_{ij}=\frac{1}{2}\wh_{ij}$, 
$\wu_{i}\cdot\ws_{ij}=\frac{1}{2}\ws_{ij}, $
$\wu_{ij}\cdot\wh_{il}=\frac{1}{2}\wh_{jl}, $
$\wu_{ij}\cdot\ws_{jl}=\frac{1}{2}\ws_{il}, $
$\wu_{ij}\cdot\wh_{ij}=\wh_{j},  $ 
$\wu_{ij}\cdot\wh_{i}=\frac{1}{2}\wh_{ij}$,
$\wh_{ij}\cdot\ws_{ij}=\frac{1}{2}(\wu_j-\wu_i)$, 
$\wh_j\cdot\ws_{jl}=\frac{1}{2}\wu_{lj}$ 
and $\wh_{ji}\cdot\ws_{il}=\frac{1}{2}\wu_{lj}$.
\end{lemma}
\proof 
The proof of the first three equalities is similar to the proof given in 
Lemma \ref{lemauihihijsij}. 
Let 
$\wu_{ij}\cdot\wh_i=\frac{1}{2}\wh_{ji}+\beta_{iji}^{ij}w_{ij}+\beta_{iji}^{ji}w_{ji}$,
$\wu_{ij}\cdot\wh_{ij}=\wh_{j}+\beta_{ij}w_j$ 
and 
$\wu_{ij}\cdot\wh_{il}=\frac{1}{2}\wh_{jl}+\beta_{ijil}^{lj}w_{lj}+\beta_{ijil}^{jl}w_{jl}$.
By the fact that $((\wu_{ij}\cdot\wh_i)\cdot\wh_i)\cdot\wu_{jl}=0$, it follows that
$\beta_{ijji}^{ij}=0$. 
Analogously  to the proof of Lemma \ref{le: regop01}, 
we get $\frac{1}{2}\beta_{ij}+\beta_{iji}^{ji}=0$.
Due to 
$((\wu_{ij}\cdot\wh_i)\cdot\wh_j)\cdot\wu_{jl}=(\wu_{ij}\cdot\wh_i)\cdot(\wu_{jl}\cdot\wh_{j})$,
it is easy to check that 
$\beta_{iji}^{ji}=0$. 
This implies that $\beta_{ij}=0$. 
Now, by \eqref{id: uijhiluiluli} 
we conclude that $\beta_{ijij}^{jl}=\beta_{ijil}^{lj}=0$.
Analogously as  Lemma \ref{le: regop01} we prove that 
$\wu_{ij}\cdot\ws_{jl}=\frac{1}{2}\ws_{il}$, 
$\wh_{ij}\cdot\ws_{ij}=\frac{1}{2}(\wu_j-\wu_i)$, 
$\wh_j\cdot\ws_{jl}=\frac{1}{2}\wu_{lj}$ 
and 
$\wh_{ji}\cdot\ws_{il}=\frac{1}{2}\wu_{lj}$, 
 and lemma follows.
\endproof

Finally, we give the proof of Theorem \ref{teo: red2}.
\proof 
By Lemmas
 \ref{le: regop01} to \ref{le: skewop01}, 
observe that if 
$\rad$ is isomorphic to one of the following radicals
$(\reg\jp)^{\mathrm{op}}$ or 
$\mathcal{P}_n$
or $(\mathcal{P}_n)^{\mathrm{op}}$
then the WPT is obvious.

Now, we assume that $\rad\cong\reg\jp$.
Let
$\wh_i$, $\wh_{ij}$, and $\ws_{ij}$ preimages of 
$\bh_i$, $\bh_{ij}$ and $\bs_{ij}$, 
respectively. 
By Lemma \ref{lemappal},
we get that 
there exist 
$\xi_{ij}$, for $i,j=1,\ldots,n$, 
such that 
\eqref{J01reg} and \eqref{J11reg} hold.

Let $\theta_1\in\corf$. 
We can choose 
$\theta_2$, $\theta_3,\ldots$, $\theta_n\in\corf$ 
such that 
$\theta_2-\theta_1=\xi_{12}-\xi_{21}$, 
$\theta_3-\theta_2=\xi_{23}-\xi_{32}$, inductively, 
$\theta_{i+1}-\theta_i=\xi_{ii+1}-\xi_{i+1i} $.
Using 
$\displaystyle{\beta_{ijjl}=\eta_{ijil}=\frac{1}{2}(\eta_{il}-\eta_{jl})}$ and 
\eqref{ec:blijbjlibijl}, it is easy to check that  
$\theta_i-\theta_j=\xi_{ji}-\xi_{ij}$ for all $i\neq j$. 
 
Let $\hh_i=\wh_i+\theta_ig_i$, 
$\hh_{ij}=\wh_{ij}+(\theta_j-\xi_{ij})g_{ij}$ and 
$\hs_{ij}=\ws_{ij}+(\xi_{ij}-\theta_j)z_{ij}\in\alga_1$.
Note that $\theta_i-\xi_{ji}=\theta_j-\xi_{ij}$ implies that
$\hh_{ij}=\hh_{ji}$ and  $\hs_{ij}=-\hs_{ji}$.
An easy computation shows that 
$\wu_i\cdot\hh_i=\hh_i$, $\wu_i\cdot\hh_{ij}=\frac{1}{2}\hh_{ij}$,
$\wu_i\cdot\hs_{ij}=\frac{1}{2}\hs_{ij}$,
$\wu_{ij}\cdot\hh_{ij}=\hh_{j}$, $\wu_{ji}\cdot\hh_{ij}=\hh_{i}$,
$\wu_{ij}\cdot\hh_{i}=\frac{1}{2}\hh_{ij}$, 
$\wu_{ij}\cdot\hh_{il}=\frac{1}{2}\hh_{jl}$,
$\wu_{ij}\cdot\hs_{jl}=\frac{1}{2}\hs_{il}$, 
$\hh_{ij}\cdot\hs_{ij}=\frac{1}{2}(\wu_j-\wu_i)$,
$\hh_{i}\cdot\hs_{ij}=\frac{1}{2}\wu_{ji}$ and 
$\hh_{ij}\cdot\hs_{jl}=\frac{1}{2}\wu_{li}$.

Further, considering 
$\algs_{\bar{1}}=\vect\langle \hh_i,\,\hh_{ij},\,\hs_{ij}\rangle$ for all $i\neq j$, then
we obtain 
$\algs_{\bar{1}}\cong(\jp)_{\bar{1}}$. 
Now,
if we take 
$\algs=\algs_{\bar{0}}\oplus\algs_{\bar{1}}\subseteq\jor$,
 then $\algs\cong\jp$. Therefore,  
an analogue to the  WPT holds, which completes the proof. 
\endproof

\medskip
\emph{Acknowledgements.} 
The authors are grateful to the referee who carefully read the paper and made valuable suggestions.

\emph{Funding}
This work is supported by Universidad de Antioquia, CODI 2016-12949.
\bibliographystyle{amsplain}

\begin{thebibliography}{99}
\bibitem{Mol} 
Th. Molien,
On systems of higher complex numbers
(Ueber Systeme h\"{o}herer complexer Zahlen),
\textit{Ann. Math.} {\bf XLI} (1893) 83--156.

\bibitem{Mol2}
Th. Molien,
Correction to the article
``On systems of higher complex numbers''
(Berichtigung zu dem Aufsatze
``Ueber Systeme h\"{o}herer complexer Zahlen''),
\textit{Ann. Math.} {\bf XLII} (1893) 308--312.

\bibitem{Wed1} S. Epsteen, J. H. Maclagan-Wedderburn,
On the structure of hypercomplex number systems,
\textit{Trans. Am. Math. Soc.}
{\bf 6} (1905) 172--178.

\bibitem{Alb1} A. A. Albert, 
The Wedderburn principal theorem for Jordan algebras;
\textit{Ann.  Math.} {\bf{48}} (1) (1947). 1-7.

\bibitem{Pen} A. J. Penico, 
The Wedderburn principal theorem for Jordan algebras;
\textit{Trans. Am. Math. Soc.} {\bf{70}}  (1951). 404-420.

\bibitem{Azk} V. G. Ashkinuze,
\emph{A theorem on the splittability of $J$-algebras, }
\textit{Ukrain. Mat. Z. 3} (1951), 391-398.

\bibitem{Sch} R. D. Schafer
The Wedderburn principal theorem for alternative algebras;
\textit{Bull. Am. Math. Soc.} {\bf{55}}  (1949). 604-614..


\bibitem{ZelShesalternativas}
E. I. Zelmanov and I. P. Shestakov
\emph{Prime alternative superalgebras and nilpotency of the radical of a free alternative algebra}
\textit{Math. USSR-Izvestiya} {\bf 37} (1) (1991) 19-36.


\bibitem{Kac} V. G. Kac, 
\emph{Classification of simple $\mathbb{Z}$-graded Lie
	superalgebras and simple Jordan superalgebras}; 
\textit{Comm. Algebra} {\bf 5} (13) (1977)1375-1400.

\bibitem{Kap1} I. Kaplansky, 
\emph{Graded Jordan algebras}, preprint.

\bibitem{Pis}
N. A. Pisarenko,
\emph{The Wedderburn decomposition in finite dimensional alternative superalgebras,}
\textit{Algebra Logic}
{\bf 32} (4) (1993) 231--238;
translation from
\textit{Algebra Logika}
{\bf 32} (4) (1993) 428--440.

\bibitem{LopD} M. C. López Díaz, 
\emph{The {W}edderburn decomposition in finite-dimensional
              alternative superalgebras of characteristic 3,}
 \textit{Comm Algebra} {\bf 28} (9) (2000) 4211- 4218.
\bibitem{Fgg18}
F.A.  Gomez-Gonzalez, 
\emph{Wedderburn principal theorem for Jordan superalgebras I}; 
\textit{Journal of Algebra} {\bf 505} (2018), 1-32.

\bibitem{Fgg16} F. A. Gómez-González, 
\emph{Jordan superalgebras of type $\suMnm$ and the Wedderburn principal theorem.}; 
\textit{Comm Algebra} {\bf 44} (7) (2016), 2867-2886.

\bibitem{FggVr} F. A. Gómez-González and Velásquez R. 
\emph{Wedderburn principal theorem for Jordan superalgebra of type $\josnm$}. 
\textit{Algebra and discrete Mathematics}, {\bf 26} (1) (2018), 19-33.

\bibitem{zelmar2} 
C. Martinez and E. I. Zelmanov 
\emph{Representation theory of Jordan superalgebras I}; 
\textit {Trans. Amer. Math. Soc}. {\bf 362} (2)(2010) 815-846.


\bibitem{Zelss} 
E. Zelmanov,
\emph{Semisimple finite dimensional Jordan superalgebras,}
Fong, Yuen (ed.) et al., Lie algebras, rings and related topics.
Papers of the 2nd Tainan--Moscow international algebra workshop 97,
Taiwan, January 11--17, 1997.
Hong Kong: Springer (ISBN 962--430--110--7/pbk). 227--243 (2000).


\bibitem{Kan} I. L. Kantor, 
\emph{Jordan and Lie superalgebras determined by a Poisson algebra}, 
\textit{Amer. Math. Soc. Transl}. {\bf 2} (1992) 151.


\bibitem{MarSheZel}
C. Martinez, I. Shestakov and E. Zelmanov, 
\emph{Jordan superbimodules over the superalgebra $P(n)$ and $Q(n)$};
\textit{Trans. Amer. Math. Soc.} {\bf 362} (4) (2010) 2037-2051.
%
%




\end{thebibliography}

\end{document}